\documentclass[11pt,twoside]{article}
\usepackage{fancyhead}
\usepackage{psfrag}
\usepackage{epsfig}
\usepackage{graphicx}
\usepackage{amsmath}
\usepackage{amssymb}
\usepackage{amsthm}
\usepackage{subfigure}
\usepackage{cite}
\usepackage[latin1]{inputenc}

\newtheorem{proposition}{Proposition}
\newtheorem{lemma}{Lemma}
\newtheorem{theorem}{Theorem}

\begin{document}

\vspace*{5mm}

\noindent
\textbf{\LARGE On completion of latin hypercuboids \\ of order 4
\footnote{This research is partially supported by  RFBR (grants
10-01-0616-a, 10-01-0424-a)} } \thispagestyle{fancyplain}
\setlength\partopsep {0pt} \flushbottom
\date{}

\vspace*{5mm}
\noindent
\textsc{V. N. Potapov} \hfill \texttt{vpotapov@math.ncs.ru} \\
{\small Sobolev Institute of Mathematics, Novosibirsk State
University, Russia}

\medskip

\begin{center}
\parbox{11,8cm}{\footnotesize
\textbf{Abstract.} A latin hypercuboid of order $N$ is an
$N\times\cdots\times N\times k$ array  filled with symbols
from the set $\{0,\dots,N-1\} $
in such a way that  every symbol
occurs at most once in every line.
If $k=N$, such an array is a latin hypercube.
 We prove that any latin hypercuboid of order $4$ is completable
to a latin hypercube.  }
\end{center}

\baselineskip=0.9\normalbaselineskip

\section{Introduction}

 An $n$-ary operation $q:\Sigma^n\to \Sigma$,
where $\Sigma$ is a nonempty set, is called an \emph{$n$-ary
qua\-si\-group} or \emph{$n$-qua\-si\-group}
(\emph{of order $|\Sigma|$}) if 
in the equality $x_{0}=q(x_1, \ldots , x_n)$ knowledge of any $n$
elements of $x_0$, $x_1$, \ldots , $x_n$ uniquely specifies the
remaining one \cite{Bel}. The \emph{Hamming distance}
$d(\overline{x},\overline{y})$ between two words $\overline{x},
\overline{y}\in \Sigma^n$ is the number of positions in which they
differ. A function $f:\Omega \subset \Sigma^n \rightarrow \Sigma$ is
called \emph{partial $n$-ary quasigroups} if $f(\overline{x})\neq
f(\overline{y})$ for any $\overline{x},\overline{y}\in \Omega$ such
that $d(\overline{x},\overline{y})=1$.

 We deal only
with {partial $n$-ary quasigroups} $f:\Sigma^{n-1}\times \Sigma'\to
\Sigma$ where $\Sigma'\subset \Sigma$. A partial $n$-ary quasigroups
$f$ is \emph{extendable} if $f= q|_{\Sigma^{n-1}\times \Sigma'}$
where $q$ is an $n$-ary quasigroup.
 The table of an $n$-ary
quasigroup values is called a \emph{latin hypercube} (the
$n$-dimensional generalization of latin square) and the table of a
partial $n$-ary quasigroup values is called a \emph{latin
hypercuboid}. By the definitions a latin hypercuboid can be
completed to a latin hypercube if and only if this partial $n$-ary
quasigroup is extendable.

  A subset $M\subset\Sigma^n$ is called an \emph{MDS code
(with distance $2$)} if $|M|=|\Sigma|^{n-1}$ and
$d(\overline{x},\overline{y})\geq 2$ for any distinct $\overline{x},
\overline{y}\in M$.  It is clear that the graph of $n$-ary
quasigroup $f$ is an MDS code; it is denoted by $M\langle f
\rangle$. Moreover there exists one-to-one correspondence between
$n$-ary quasigroups and MDS codes.

A subset $M\subset\Sigma^n$ is called a \emph{$k$-fold MDS code } if
 for any $i\in \{1,\dots,n\}$ and word $\overline{x}\in \Sigma^n$ there
exist exactly $k$ words $\overline{y}\in M$ such that the words
$\overline{x}$ and $\overline{y}$ coincide in all non-$i$th
coordinates.  The complement $\Sigma^{n}\setminus M$ of a $k$-fold
MDS code is a $(|\Sigma|-k)$-fold MDS code. Any partial $n$-ary
quasigroup $q:\Sigma^{n-1}\times \Sigma'\to \Sigma$ corresponds to a
$k$-fold MDS code $M\subset\Sigma^{n}$, where $k=|\Sigma'|$. A
$k$-fold MDS code $M\subset\Sigma^{n}$ is called \emph{splittable}
if it can be represented as the union of $k$ disjoint MDS codes,
i.\,e. it corresponds to a partial $n$-ary quasigroup. Besides the
complement $\Sigma^{n}\setminus M$ is splittable if and only if this
partial $n$-ary quasigroup is extendable.

  The fact that every
latin rectangle can be completed to a latin square is a simple
consequence of Konig's theorem (see \cite{Hol}). It is well-known
(see \cite{Kochol}) that any partial $n$-ary quasigroup
$f:\Sigma^{n-1}\times \Sigma'\to \Sigma$ is extendable if
$|\Sigma'|=1$ or $|\Sigma'|=|\Sigma|-1$. Therefore all partial
$n$-ary quasigroups of order $N=|\Sigma|\leq 3$ are extendable.

Kochol (\cite{Kochol} and \cite{Kochol1}) proved that for any $k$
and $N$ satisfying $N/2<k<N-2$ there is an $N \times N\times k$
latin cuboid that cannot be completed to a latin cube. Examples of
non-completable $N \times N\times k$ latin cuboids  for $N=5,6,7,8$
and $k=2,2,3,4$ respectively are constructed  in
 \cite{MKW}. Using these facts it is easy to show that there are $n$-dimensional latin
hypercuboids of order $N$ that cannot be completed to a latin
hypercube for $n\geq 3$ and $N\geq 5$. In \cite{KP04} we constructed
examples of unsplittable $2$-fold MDS codes for $n\geq 3$ and $N=4$.
But the question of existence of a non-completable latin cuboid for
$N=4$ was open.

\section{Main result }

In this paper we announce the following theorem.

 \begin{theorem} Any latin hypercuboid of order $4$ is completable to a latin
  hypercube. \end{theorem}

 For the proof of the theorem  we need additional
  notations.

  Sets $M_1, M_2\subset \Sigma^n$ are called {\em isotopic} if for
some collection of $n+1$ permutations ({\em isotopies}) $\bar
\tau=(\tau_0,\tau_1,\ldots ,\tau_n)$,
$\tau_i:\Sigma\rightarrow\Sigma$, we have
$\chi_{M_1}(x_1,\dots,x_{n}) \equiv \chi_{M_2}(\tau_1 x_{1}, \ldots,
\tau_{n} x_{n})$, where $\chi_{M}$ is an indicator of the set $M$.

$n$-Ary quasigroups $f$ and $g$ are called {\em isotopic} if MDS
codes $M\langle f\rangle$ and $M\langle g\rangle$ are isotopic.

$n$-Ary quasigroups $f$ and $g$ are called {\em parastrophic}  if
$$
 x_{\sigma0}=f(x_{\sigma1},\ldots,x_{\sigma n})\ \Leftrightarrow \  x_{0}=g(x_1,\ldots,x_n)
$$
for some  permutation $\sigma: \{1,\ldots,n\} \to \{1,\ldots,n\}$.

 If we assign some fixed values to $l\in \{1,\ldots,n-1\}$ variables
 of an
 $n$-ary quasigroup $f$ or any of its parastrophs,  then we obtain
an $(n-l)$-qua\-si\-group. Such  qua\-si\-groups are called
\emph{retracts} or \emph{$(n-l)$-retracts} of $q$. If $x_0$ is not
fixed, the retract is \emph{principal}.

 An $n$-qua\-si\-group $f$ is termed {\em permutably reducible}
if there exist \linebreak $m\in\{2,\ldots,n-1\}$, an
$(n-m+1)$-qua\-si\-group $h$, an $m$-qua\-si\-group  $g$, and a
permutation $\sigma: \{1,\ldots,n\} \to \{1,\ldots,n\}$ such that
$$f(x_1,\ldots,x_{n}) \equiv h(g(x_{\sigma(1)},\ldots, x_{\sigma(m)}),
 x_{\sigma(m+1)},\ldots, x_{\sigma(n)}).$$
For short, we will omit the word ``permutably''. If an
$n$-qua\-si\-group is not reducible, then it is {\em irreducible}.
By the definition, all binary quasigroups are irreducible.

An $n$-qua\-si\-group $f$ is termed {\em completely reducible} if it
is reducible and all its principal retracts of arity more than $2$
are reducible.

Let $\Sigma=\{0,1,2,3\}$. Consider the $2$-fold MDS code $L\subset
\Sigma^n$ defined by indicator function
$\chi_L(x_1,\ldots,x_n)\equiv
\chi_{0,1}(x_1)\oplus\cdots\oplus\chi_{0,1}(x_n)$. We say that a
$2$-MDS code is {\em linear} if it is isotopic to 2-MDS code $L$.
 Let $f$ be an $n$-ary  quasigroup and $a,b\in \Sigma, a\neq b$. Then the set
 $S_{a,b}(f)=\{(x_1,\ldots,x_n)\in \Sigma^n \ | \
f(x_1,\ldots,x_n)\in \{a,b\}\}$ is a $2$-fold MDS code. An $n$-ary
quasigroup $f$ is called {\em semilinear} if for some $a,b\in
\Sigma$ the set $S_{a,b}(f)$  is linear.

The proof of the theorem  is based on the following statements.

\begin{proposition}{\rm \cite{KP09}}
 Every
$n$-ary quasigroups of order $4$ is reducible or  semilinear.
 \end{proposition}

\begin{proposition}{\rm \cite{Cher}} Any reducible $n$-ary quasigroup $f$  is
representable in the form
\begin{equation}
 f(\overline{x})\equiv q_0(q_2(\tilde x_1),...,q_m(\tilde
x_m)),\label{eq:Decomp.-of-quas.}
\end{equation}
 where $q_j$ are $n_j$-ary quasigroups where $1\leq j\leq m$, $q_0$
 is an
irreducible an $m$-ary quasigroup,
 $ \tilde x_j$ ---  disjoint collections of variables from the set
 $\{x_i\}$. If $m\geq 3$ there is only one collection of the sets $ \tilde
 x_j$ meeting equation {\rm(\ref{eq:Decomp.-of-quas.})}.
 \end{proposition}

\begin{proposition}{\rm \cite{Cher}} Any completely reducible $n$-ary quasigroup $f$ of order
$4$ is representable in the form
 \begin{equation}
 f(\overline{x})\equiv q_1(\tilde x_1)\ast...\ast
q_k(\tilde x_k),\label{eq:Decomp.-of-quas1.}
\end{equation}
where $\ast$ is a group operation, $q_j$ are $n_j$-ary quasigroups
$(1\leq j\leq k)$ that cannot be decomposed as $q_j(\tilde x_j)
\equiv q'(\tilde x'_j)\ast q''(\tilde x''_j)$,
 $ \tilde x_j$ --- disjoint collections of variables from the set $\{x_i\}$.
 There is only one collection of the sets $ \tilde x_j$ meeting
 equation
  {\rm(\ref{eq:Decomp.-of-quas1.})}.
\end{proposition}

We say that decomposition of $n$-ary quasigroup $f$ represented
above is {\em canonical}.

As noted above, it is sufficient to prove that any two disjoint
  $n$-ary quasigroups of order $4$ can be completed to an $(n+1)$-ary
  quasigroup.
  We will argue this statement by induction on $n$. For $n=1,2,3,4$
  this fact can be verifed by computer.
  Let the induction hypothesis be true for every integer $m$, $m\leq n$.
 Induction step is proved by the following lemmas.

\begin{lemma}
  Let $f$ and  $g$  be   $n$-ary quasigroups, let
$f(x_1,\ldots,x_n)\neq g(x_1,\ldots,x_n)$ for all
$(x_1,\ldots,x_n)\in \Sigma^n$, and  let one (or both) of these
quasigroups be irreducible. Then there exists an $(n+1)$-ary
quasigroup $F$ such that $F|_{x_{n+1}=0}= f$ and $F|_{x_{n+1}=0}=
g$.
\end{lemma}

\begin{lemma} Let $f$ and  $g$  be reducible  $n$-ary quasigroups (but one of
them is not completely reducible), $f(x_1,\ldots,x_n)\neq
g(x_1,\ldots,x_n)$ for all $(x_1,\ldots,x_n)\in \Sigma^n$. Let the
canonical decompositions of $f$ and $g$ be
$$ f  (\overline{x})\equiv
q_0(q_1(\tilde x_1),...,q_m(\tilde x_m))\}, $$
$$g(\overline{x}) \equiv
q'_0(q'_1(\tilde x'_1),...,q'_{m'}(\tilde x'_{m'}))\},$$ where
disjoint collections of variables $\{\tilde x_i\}$ and  $\{\tilde
x'_i\}$ are different. Then there exist an $(n+1)$-ary quasigroup
$F$ such that $F|_{x_{n+1}=0}= f$ and $F|_{x_{n+1}=0}= g$.
\end{lemma}

 \begin{lemma} Let $f$ and  $g$  be reducible  $n$-ary quasigroups,
$f(x_1,\ldots,x_n)\neq g(x_1,\ldots,x_n)$ for all
$(x_1,\ldots,x_n)\in \Sigma^n$. Let
$$M\langle f \rangle= \{x\in \Sigma^{n+1} :  q_1(\tilde x_1)=
q_0(q_2(\tilde x_2),...,q_m(\tilde x_m))\}, $$
$$M\langle g \rangle= \{x\in \Sigma^{n+1} :  q'_1(\tilde x_1) =
q'_0(q'_2(\tilde x_2),...,q'_{m}(\tilde x_{m}))\},$$ where $q_0$ and
$q'_0$ are irreducible $m$-ary quasigroups, $m\geq 4$,  $\tilde x_i$
--- disjoint collections of variables. Then there exists an
$(n+1)$-ary quasigroup $F$ such that $F|_{x_{n+1}=0}= f$ and
$F|_{x_{n+1}=0}= g$.
\end{lemma}

  \begin{lemma}
  Let $f$ and  $g$   be  completely reducible $n$-ary quasigroups,
$f(x_1,\ldots,x_n)\neq g(x_1,\ldots,x_n)$ for all
$(x_1,\ldots,x_n)\in \Sigma^n$.  Then there exists an $(n+1)$-ary
quasigroup $F$ such that $F|_{x_{n+1}=0}=f$ and $F|_{x_{n+1}=0}= g$.
 \end{lemma}

 It is absolutely clear that we have considered all the cases.

\end{document}